\documentclass[12pt]{article}
\usepackage{amsfonts,url}
\usepackage{amsthm,latexsym,amssymb}
\pdfoutput=1
\usepackage{amsmath}
\usepackage{graphicx}
\usepackage{comment}
\usepackage[usenames,dvipsnames]{xcolor}
\usepackage{enumerate}
\usepackage{extarrows,pdfsync}
\usepackage{hyperref}
\usepackage[sort]{natbib}
\usepackage{stackrel}
\usepackage{centernot}
\usepackage{caption}
\usepackage{subcaption}
\usepackage{fancyhdr}
\usepackage{geometry}
\let\quoteOLD\quote 
\def\quote{\quoteOLD\small}

\definecolor{labelkey}{cmyk}{0,0.8,1,0.5}
\definecolor{refkey}{cmyk}{0,0.8,1,0.5}

\newtheorem{theorem}{Theorem}[section]
\newtheorem{example0}{\sc Example}[subsection]

\newtheorem{eg}{Example}

\newtheorem{corollary}{Corollary}
\newtheorem{proposition}{Proposition}

\newtheorem{remark}{Remark}

\numberwithin{equation}{section}
\numberwithin{theorem}{section}
\numberwithin{corollary}{section}
\numberwithin{proposition}{section}
\numberwithin{lemma}{section}
\numberwithin{definition}{section}
\numberwithin{remark}{section}

\makeatletter
\def\th@newremark{\th@remark\thm@headfont{\bfseries}}
\makeatletter

\def\boxit#1{\vbox{\hrule\hbox{\vrule\kern6pt
          \vbox{\kern6pt#1\kern6pt}\kern6pt\vrule}\hrule}}

\newcommand{\pibar}{\overline{\Pi}}
\newcommand{\pibarinv}{\overline{\Pi}^{\leftarrow}}

\newcommand{\lchinv}{{h^\leftarrow}}
\newcommand{\hinv}{{H^\leftarrow}}

\newcommand{\beqq}{\begin{equation}}
\newcommand{\eeqq}{\end{equation}}

\newcommand{\eqdr}{\stackrel{\mathrm{D}}{=}}

\newcommand{\R}{\Bbb{R}}

\newcommand{\N}{\Bbb{N}}

\newcommand{\rmd}{{\rm d}}
\newcommand{\rmi}{{\rm i}}
\newcommand{\halmos}{\quad\hfill\mbox{$\Box$}}

\newcommand{\XX}{\mathbb{X}}

\newcommand{\dto}{\downarrow}

\newcommand{\be}{\begin{equation}}
\newcommand{\ee}{\end{equation}}
\newcommand{\bea}{\begin{eqnarray}}
\newcommand{\eea}{\end{eqnarray}}
\newcommand{\bean}{\begin{eqnarray*}}
\newcommand{\eean}{\end{eqnarray*}}
\newcommand{\ben}{\begin{equation*}}
\newcommand{\een}{\end{equation*}}
\newcommand{\ba}{\begin{aligned}}
\newcommand{\ea}{\end{aligned}}

\def\nexto{\kern -0.54em}

\newcommand{\PP}{\textbf{\rm P}}
\newcommand{\EE}{\textbf{\rm E}}

\newcommand{\bm}[1]{\mbox{\boldmath $#1$\unboldmath}}

\newcommand{\BGamma}{\bm \Gamma}
\newcommand{\BPi}{\bm \Pi}

\newcommand{\nuu}{\Pi}
\newcommand{\Q}{\pibar}
\newcommand{\Y}{\Delta^{(r)}}
\newcommand{\trimX}{{}^{(r)}X}

\begin{document}

\title{\bf Trimmed L\'evy Processes and their Extremal Components}
\author{Yuguang Ipsen, Ross Maller and Sidney Resnick\thanks{Research
partially supported by ARC Grants DP1092502 and DP160104737. S. Resnick was also partly
supported by Army MURI grant 
  W911NF-12-1-0385 to Cornell University and acknowledges with thanks
  hospitality and support from the FIM - Institute for Mathematical
  Research, ETH Zurich in April 2017. Resnick also acknowledges
  hospitality from the Department 
  of Mathematics, Oregon State University in January 2018.
  \newline
Email: Yuguang.Ipsen@anu.edu.au; Ross.Maller@anu.edu.au; sir1@cornell.edu}}


\maketitle

\begin{abstract}
We analyse a trimmed stochastic process of the form 
 ${}^{(r)}X_t= X_t - \sum_{i=1}^r \Delta_t^{(i)}$,
 where $(X_t)_{t \geq 0}$ is a driftless subordinator on $\R$ 
with its jumps on $[0,t]$ ordered as
$ \Delta_t^{(1)}\ge  \Delta_t^{(2)} \cdots$.
When  $r\to\infty$, both $\trimX_t\dto 0$ and $\Y_t\dto 0$ a.s. for each $t>0$,  and it is interesting to study the  weak limiting behaviour of $\bigl(\trimX_t, \Y_t\bigr)$ in this case.
We term this ``large-trimming" behaviour.
Concentrating on the case $t=1$, we study joint convergence of
$\bigl(\trimX_1, \Y_1\bigr)$  under linear normalization,
assuming extreme value-related conditions on the L\'evy measure of $X$
which  guarantee that $\Y$ has a limit distribution with linear
normalization.  
Allowing $\trimX$ to have random centering and scaling in  a natural way, we show that $\bigl(\trimX_1, \Y_1\bigr)$  has a bivariate normal
limiting distribution, as $r\to\infty$; but 
replacing the random normalizations with natural deterministic ones produces non-normal limits which we can specify.
\end{abstract}

extreme value-related conditions, large-trimming limits.


\section{Introduction}\label{s1}
%

Suppose $(X_t)_{t \geq 0}$ is a driftless subordinator
with infinite L\'evy measure $\nuu $ 
and tail function 
$\Q(x):=\nuu(x,\infty), \, x>0$.  
Thus, $(X_t)$ has Laplace transform 
$Ee^{-\lambda X_t}= e^{-t\psi(\lambda)}$, $t\ge 0$, where 
\ben
\psi(\lambda)= \int_{(0,\infty)}(1-e^{-\lambda x})\Pi(\rmd x),\ \lambda>0.
\een

Let $\Y_t$ be the $r$th largest jump of $X_t$ on
$[0,t]$, $t>0$, $r\in\N:=\{1,2,\ldots\}$. 
The trimmed subordinator is defined to be 
 ${}^{(r)}X_t= X_t - \sum_{i=1}^r \Delta_t^{(i)}$,  $t>0$, $r\in\N$. 
 In \cite{buchmann:maller:resnick:2016, ipsen:maller:resnick:2017} we considered distributional
 properties of  
$\Delta_t^{(r)}$
 as a function of $r$ and here we continue that study by considering
 the joint weak limiting behaviour of $\bigl(\trimX_t, \Y_t\bigr)$ as
 $r\to\infty$.  
  As $r\to\infty$, $\trimX_t\dto 0$ and $\Y_t\dto 0$ a.s. for each
  $t>0$, 
but conditionally on $\Y_t$ we may consider $\trimX_t$ as a L\'evy
process with L\'evy measure restricted to $(0,\Y_t)$
(e.g., \cite{resnick:1986b}).
So as $r\to\infty$ and big jumps are removed from $\trimX$, it makes
sense that we expect $\trimX$ should have a Gaussian weak limit
after centering and norming.
We focus on the case
 $t=1$ and write simply $\bigl(\trimX, \Y\bigr)$ for $\bigl(\trimX_1,
 \Y_1\bigr)$
 (the case of general $t>0$ is considered briefly in Section \ref{s7}).

%

 The approach we take is to assume conditions on $\Q$ guaranteeing
 that $\Y$ has a limit distribution under linear normalization, and
 then prove that a normal limit distribution of $\trimX$ conditional
 on the value of $\Y$ also exists as $r\to\infty$. For finite $r$, we
 denote the conditional distribution  with the notation $\trimX | \Y$.  The
 conditioned limit of $\trimX | \Y$ initially requires a natural
 random centering and random scaling to achieve asymptotic normality.
 Having derived that, we then investigate replacing the random
 centering and scaling with deterministic versions.

According to \cite[Section 4.2]{buchmann:maller:resnick:2016},
there exist scaling functions $a_r>0$ and centering functions $b_r \in
\mathbb{R}$  such that, as $r \to \infty$, weak convergence\footnote{We use the symbol ``$ \Rightarrow$" to denote weak convergence in $\R$ or $\R^2$.}
 holds in $\R$:
\begin{equation}\label{Yconv}
\frac{\Y  -b_r}{a_r} \Rightarrow \Delta^{(\infty )},
\end{equation} 
with $\Delta^{(\infty )} $ non-degenerate, 
if,  for $x\in\R$ such that $a_rx+b_r >0$,
\begin{equation}\label{e:LTcondit}
\lim_{r\to\infty}
\frac{ r-\Q\bigl(a_rx+b_r \bigr)}{\sqrt r}=h(x),
\end{equation}
where $h(x) \in \R$ is a 
non-decreasing and non-constant limit function.
The function $h(x)$ has the form (see \cite{buchmann:maller:resnick:2016}, Eq. (4.2)):
\beqq
\label{e:h}
\tfrac 12 h(x)=\tfrac 12 h_\gamma(x)= 
\begin{cases}
-\frac{1}{\gamma} \log (1-\gamma x),& \text{ if }\gamma\in\R\setminus\{0\},\,1-\gamma x>0,\\
x,& \text{ if } \gamma=0,\,x\in\R.
\end{cases}
\eeqq

We can identify the distribution of the limit  random variable  $\Delta^{(\infty )}$
in terms of the inverse function $\lchinv$ of $h$. 
From \eqref{e:h} this function satisfies,  for $y\in \R$,
\begin{equation}\label{e:hinv}
h^\leftarrow(y)=\frac{1-e^{-\gamma y/2}}{\gamma} =
\begin{cases}
y/2,& \text{ if }\gamma=0,\\[5pt] 
\displaystyle{\frac{1-e^{-\gamma y/2}}{\gamma}},& \text{ if }\gamma>0,\\[10pt] 
\displaystyle{\frac{e^{|\gamma |y/2} -1}{|\gamma|}},& \text{ if }\gamma<0.
\end{cases}
\end{equation}
We note that $h^\leftarrow : \R \mapsto \R_\gamma$, where, for $\gamma \in \R$, 
$$\R_\gamma :=\{x\in \R: 1-\gamma x>0\}=\begin{cases}
\R,& \text{ if }\gamma=0,\\
(-\infty, \frac{1}{\gamma}),& \text{ if } \gamma>0,\\
(-\frac{1}{|\gamma |}),\infty),& \text{ if } \gamma<0. \end{cases} $$

Taking inverses in \eqref{e:LTcondit}, we get an equivalent form
\beqq\label{e:Qinv}
\lim_{r\to\infty}\frac{\Q^\leftarrow (r-y\sqrt r) -b_r}{a_r}
= h^\leftarrow (y), \quad \quad y \in \R,
\eeqq
where the inverse function $\pibarinv$ to $\pibar$ is defined by
\ben
\pibarinv(x)= \{\inf y>0:\pibar(y)\le x\}.
\een
From \eqref{e:hinv} we have $h^\leftarrow (0)=0$, so from 
\eqref{e:Qinv} we deduce for  $y \in \R,$
\bea\label{r4}
\lim_{r\to\infty}\frac{\Q^\leftarrow (r-y\sqrt r) -\Q^\leftarrow (r)}{a_r}
&=&
\lim_{r\to\infty} \Bigl(\frac{\Q^\leftarrow (r-y\sqrt r) -b_r}{a_r}
-\frac{\Q^\leftarrow (r) -b_r}{a_r}\Bigr)\cr
&&\cr
&=&
h^\leftarrow (y)- h^\leftarrow (0)= h^\leftarrow (y) .
\eea
We conclude that for centering constants we may always set 
$b_r=\Q^\leftarrow (r)$
(appropriate norming constants $a_r$ will be specified later). 

The convergences in \eqref{e:LTcondit}, \eqref{e:Qinv} and \eqref{r4} are locally uniform since they are convergences of monotone functions to a continuous limit. 
Recalling the notation in \eqref{Yconv}, we have under
\eqref{e:LTcondit} that
\begin{equation}\label{e:onedim}
\lim_{r\to\infty} \PP\Bigl( \frac{\Y  -b_r}{a_r} \leq x \Bigr) = \PP\bigl
( \Delta^{(\infty )} \leq x \bigr)=\Phi\bigl( h(x)\bigr),\ x\in\R,
\end{equation}
where $\Phi(x)$ is the standard normal cdf. (This will be proved in  \eqref{e:YAN} below, or see \cite{buchmann:maller:resnick:2016}). Thus $\Delta^{(\infty )}
\eqdr  h^\leftarrow (N(0,1))$ where $N(0,1)$ is a standard normal random variable. 

\begin{remark}\label{rem:noPosGamma}
{\rm 
Since we assume only positive jumps for the L\'evy process,
$\Pi(\cdot)$ concentrates on $(0,\infty)$. This implies that the case
$\gamma>0$ in \eqref{e:h} or \eqref{e:hinv} cannot occur. 
From the discussion in
\cite{buchmann:maller:resnick:2016}, \eqref{e:LTcondit} means that the function  $G(x):=e^{-\sqrt{\Q(x)}}$ defined on $(0,\infty)$ is a distribution function in the minimal domain
 of attraction, which for the $\gamma>0$ case would require $G(x)$ to be regularly varying as  $x\to -\infty$. This is impossible because $\Pi(\cdot)$ concentrates on $\R_+$.
 So from now on we concentrate attention on the cases $\gamma\le 0$.
 }
\end{remark}

To conclude this introduction we set out the steps we intend to follow to understand the joint limit behaviour of
$(\trimX, \Y)$ as $r\to\infty$ under \eqref{e:LTcondit} or, equivalently,  \eqref{e:Qinv}.

\begin{enumerate}
\item As discussed, we expect a normal limit as $r\to\infty$ for $\trimX$ with   suitable linear normalizations. We show that this
  happens for $\trimX |\Y$ under a natural {\it random} centering and scaling (Theorem \ref{thm:trimXAN}).
\item Following that, we extend asymptotic normality of  $\trimX |\Y$ 
to a joint asymptotic weak limit for $(\trimX, \Y)$ in which the limit has independent components. At this stage,   $\trimX$ still has random centering and scaling, though $\Y$ has non-random normalizations
(Corollary \ref{cor:joint}).
\item Finally, we note there is a cost to replacing the random centering and scaling:  dependencies and non-normality  are introduced into the limit  (Theorems \ref{thm:joint} and \ref{e:repmu}).
\end{enumerate}

In the next section we give our main results. Proofs of the theorems and
further discussion are deferred to Sections \ref{sect:proofs} and \ref{p23}. A number of subsidiary propositions are also needed; these are proved in Sections \ref{s3} and \ref{sec:more}. Section \ref{s7} concludes with some general discussion.

%
\section{Main Results}\label{s2}
Throughout, we write $\PP^{\Y} (\cdot) = \PP( \cdot |\Y )$
for  the conditional distribution, given $\Y $.
In introducing this we make the simplifying assumption that $\Pi$ is atomless (equivalently, $\pibar$ is continuous on $(0,\infty)$).
This means that the inverse function $\pibarinv$ is strictly increasing on $(0,\infty)$ and the ordered jumps $\Delta X_t^{(i)}$ are uniquely defined.
We expect that this assumption can be removed by some well known manipulations which would add little of interest to the exposition, so we omit them.

We will also need  truncated first and second moment functions,
 defined for $t>0$ by 
\beqq \label{e:defmusigma}
\mu(t)=\int_0^{t} x\, \nuu (\rmd x)\quad \text{and} \quad 
\sigma^2(t)= \int_0^{t} x^2 \, \nuu (\rmd x).
\eeqq

\begin{theorem} \label{thm:trimXAN}
Suppose $X$ is a 
driftless subordinator on   $(0,\infty)$ with L\'evy measure $
\nuu(\cdot)$ on $(0,\infty)$ that satisfies \eqref{e:LTcondit} or, equivalently, \eqref{e:Qinv}, 
 for deterministic functions $a_r > 0$ and $b_r \in \R$. 
Then 
we have
\beqq
\lim_{r\to\infty}
\PP^{\Y} \Bigl( \frac{ \trimX -\mu (\Y)}{\sigma (\Y)} \leq x\Bigr) =
\Phi (x), \ x\in\R.\label{e:XCONV}
\eeqq
\end{theorem}

\begin{remark}\label{r2}
{\rm 
By the dominated convergence theorem the convergence in \eqref{e:XCONV} holds unconditionally as well, so we also have
\ben
 \frac{ \trimX -\mu (\Y)}{\sigma (\Y)} \Rightarrow
N(0,1),\ {\rm as}\ r\to\infty,
\een
under the conditions of Theorem \ref{thm:trimXAN}.
}
\end{remark}

Retaining the random centering and scaling, 
Theorem \ref{thm:trimXAN} immediately leads to a joint limit
distribution for $(\trimX,\Y)$.  
In the  following corollary,  $N_X$ and $N_\Gamma$ are independent standard normal random variables, being the limits of the standardised $\trimX$ and $\Y$, with the subscripts on $N_X$ and $N_\Gamma$ serving to distinguish the components  corresponding to $\trimX$ and $\Y$.
(Throughout, $N_X$ and $N_\Gamma$ will be independent standard normal random variables corresponding to $\trimX$ and $\Y$ in this way.)

\begin{corollary}\label{cor:joint} Under the conditions leading to
  \eqref{e:onedim} and \eqref{e:XCONV} we have, in $\R^2$,
$$\Bigl(
\frac{\trimX-\mu(\Y)}{\sigma(\Y)}, \,
\frac{\Y-b_r}{a_r} \Bigr) \Rightarrow 
\bigl(N_X,\,  h^\leftarrow(N_\Gamma)\bigr),\ {\rm as}\ r\to\infty.
$$
\end{corollary}

Next we need to understand the effect of replacing the random
centering and scaling by deterministic counterparts.
We begin with the scaling constants. The treatment is broken up according to the
cases of the constant $\gamma $ in \eqref{e:h}.

\begin{theorem}\label{thm:joint}
Suppose \eqref{e:LTcondit} holds.

\begin{enumerate}[\rm (i)]
\item When $\gamma < 0$,  we have, as $r\to\infty$, with $b_r=\Q^\leftarrow (r)$,
\be\label{joint:<0}
\Bigl(
\frac{\trimX-\mu(\Y)}{\sigma(b_r)}, \,
\frac{\Y}{b_r} \Bigr)
\Rightarrow 
\Bigl(N_X e^{-N_\Gamma |\gamma|/2},\, e^{-N_\Gamma |\gamma|/2}\Bigr)
\ee
and removing the random centering from $\trimX$ gives
\be\label{r1}
\Bigl(
\frac{\trimX-\mu(b_r)}{b_r\sqrt r},\, 
\frac{\Y}{b_r} \Bigr)
\Rightarrow 
\Bigl(\frac{2}{|\gamma|}(e^{-N_\Gamma |\gamma|/2} -1),\,
 e^{-N_\Gamma |\gamma|/2}\Bigr),\ {\rm as}\ r\to\infty.
\ee

\item 
When $\gamma=0$,   we have, as $r\to\infty$, with  $a_r = 2(\pibarinv(r-\sqrt{r}) - \pibarinv(r))$ and $b_r = \pibarinv(r)$, 
\be\label{e:gamma0}
\Bigl(
\frac{\trimX-\mu(\Y)}{\sigma(b_r)}, \,
\frac{\Y-b_r}{a_r} \Bigr) \Rightarrow \Bigl(N_X,\, \frac{N_\Gamma}{2}\Bigr),\ {\rm as}\ r\to\infty.
\ee
\end{enumerate}
\end{theorem}

\begin{remark}
{\rm
\begin{enumerate}[(a)]
\item Note that when $\gamma<0$,  we no longer have independence of the components in the limit when we
replace the random scaling by the deterministic one as in \eqref{joint:<0} and \eqref{r1}.

\item 
When $\gamma = 0$, we can always make the scaling deterministic, as in \eqref{e:gamma0},
however this is not in general the case for the centering;
replacing $\mu(\Y)$ with $\mu(b_r)$ in \eqref{e:gamma0} is only possible under some subsidiary conditions.
A detailed discussion is given in Section \ref{p23}. For the special case when $\pibar \in RV_0(-\alpha)$ for  $0\le \alpha\le 1$, 
the joint limiting distribution of $\trimX$ and $\Y$ is specified in the following theorem. 
\end{enumerate}
}
\end{remark}

\begin{theorem}\label{e:repmu}
Suppose $\gamma=0$ and 
$\Q$ is regularly varying at $0$ with index  $-\alpha$.
Let $c_\alpha:= \alpha/(2-\alpha)$,   $a_r = 2(\pibarinv(r-\sqrt{r}) - \pibarinv(r))$ and $b_r = \pibarinv(r)$.

(i)\  Suppose  $0<c_\alpha\le 1$, so that $\alpha\le 1$.
Then
\beqq\label{e:depCase}
\Biggl(
\frac{\trimX-\mu(b_r)}{\sigma(b_r)}, \,
\frac{\Y-b_r}{a_r} \Biggr) \Rightarrow 
\Bigl(N_X +\frac{N_\Gamma}{\sqrt{c_\alpha}},\, \frac{N_\Gamma}{2}\Bigr).
\eeqq


(ii)\  Suppose  $c_\alpha= 0$, so that $\pibar$ is slowly varying at 0.
Then 
$$
\Bigl(\frac{\trimX-\mu(b_r)}{ b_r\sqrt{r}}, \,
\frac{\Y-b_r}{a_r} \Bigr) \Rightarrow 
\Bigl(N_X, \, \frac{N_\Gamma}{2}\Bigr).
$$
\end{theorem}

\section{Convergence of $\Y$}\label{s3}
We begin the program outlined in the previous section by examining the convergence of $\Y$, after norming and centering, as $r\to\infty$. 
Throughout this section assume\footnote{The simplifying assumption that $\pibar$ is continuous on $(0,\infty)$ is not in fact needed for the results in this section.} \eqref{e:LTcondit} or, equivalently,  \eqref{e:Qinv}, and recall the function $h$ in \eqref{e:h} and its inverse $\lchinv$ in \eqref{e:hinv}.

We need some more preliminary setting up.
Let $\{\Gamma_l\}$ and $\{\Gamma'_l\}$ be cumulative sums of
independent sequences of iid standard exponential random variables.
We can construct the subordinator $X$ from a Poisson random measure  
$$\XX(\cdot) =\sum_{l=1}^\infty \delta_{\Q^\leftarrow (\Gamma_l)}$$
where the mean measure is $\nuu (\cdot)$ and the points are written in decreasing order
\cite{ferguson:klass:1972}, \cite{lepage:woodroofe:zinn:1981},
\cite{resnick:1986b}, \cite{buchmann:fan:maller:2016}
and \cite[p.139, Ex. 3.38]{resnickbook:2008}.
This means
$$X=\int_0^\infty x \, \XX(\rmd x)=\sum_{l=1}^\infty \Q^\leftarrow (\Gamma_l)
\quad {\rm and}\quad 
\Y =\Q^\leftarrow (\Gamma_r);
$$
also
$$\trimX
= \sum_{l=r+1}^\infty \Q^\leftarrow (\Gamma_l) = \sum_{l=1}^\infty \Q^\leftarrow (\Gamma_r +\Gamma'_l).
$$

\subsection{Proof of the convergence in (\ref{e:onedim})} 
We may understand the form of the limit for $\Y$ in \eqref{e:onedim} as follows.   By properties of the gamma distribution,  we know that, as $r\to\infty$,  $G_r:=(\Gamma_r -r)/\sqrt r
 \Rightarrow N(0,1)$, a standard normal random variable. Assume  \eqref{e:LTcondit},  or, equivalently,  \eqref{e:Qinv};
 then, owing to the local uniform convergence in \eqref{e:Qinv},
we get \eqref{e:onedim} from 
\begin{align}
\frac{\Y-b_r}{a_r} =&
\frac{\Q^\leftarrow (\Gamma_r)-b_r}{a_r}
=\frac{\Q^\leftarrow (r+G_r \sqrt r)-b_r}{a_r}
 \Rightarrow 
   h^\leftarrow (-N(0,1))\eqdr    h^\leftarrow
   (N(0,1)).\label{e:YAN}
\end{align}


\subsection{Role of the de Haan classes $\BGamma$ and $\BPi$}
\label{subsec:implic} 

Now introduce the function $H:[0,\infty)\mapsto [1,\infty)$ defined by
\beqq\label{e:H}
H(t)=e^{2 \sqrt t}, t>0,
\eeqq
and define the non-increasing function $V$ by
\beqq\label{e:defV}
V(x)=\Q^\leftarrow \circ H^\leftarrow (x),\quad x>1, 
\eeqq
and changing variables gives the  representation $\pibarinv(x)=V(H(x))$, 

The function $H$ is the canonical example of a non-decreasing
 function in the de Haan class
${\bf \Gamma}$ with auxiliary function $f(t)=\sqrt t$ \cite{dehaan:1970,bingham:goldie:teugels:1989,
dehaan:1974,dehaan:resnick:1973,geluk:dehaan:1987,resnickbook:2008}
satisfying
\beqq\label{e:limH}
\lim_{t\to\infty}
\frac{H(t+xf(t))}{H(t)}=e^x, \quad 
x\in\R.
\eeqq
This can be verified directly or by reference to \cite[p. 248, line -1]{dehaan:1974}.
The inverse function $H^\leftarrow:[1,\infty)\mapsto
[0,\infty)$ to $H$ is $H^\leftarrow
(y)=\frac 14 \log^2 y,\,y>1$, and inverting \eqref{e:limH}  shows that
$H^\leftarrow$ satisfies
\beqq\label{e:limHinv}
\lim_{s\to\infty} \frac{H^\leftarrow (sy) -H^\leftarrow
  (s)}{f(H^\leftarrow (s))} =\log y, \quad y>0,
\eeqq
so $H^\leftarrow$  is an increasing function in de Haan's function
class ${\bf{\Pi}}$
(\cite{dehaan:1970,bingham:goldie:teugels:1989,resnickbook:2008} or
\cite[p. 375]{dehaan:ferreira:2006}).
It has slowly varying  
auxiliary function $g(s)=f\circ H^\leftarrow (s)=\sqrt{H^\leftarrow (s)}=\frac 12 \log s$ which is the denominator in
\eqref{e:limHinv}. The convergence in \eqref{e:limHinv} is uniform in compact intervals of $y$ bounded away from 0.

Recall that \eqref{e:Qinv} is in force throughout, so we have \eqref{r4} also.
Applying the uniform convergence in \eqref{r4} and
\eqref{e:limHinv},  we see that $V$ satisfies, for $x>0$,
\bean
\lim_{s\to\infty} 
\frac{V (sx) -V(s)}{a\circ H^\leftarrow (s)}
&=&
\lim_{s\to\infty} 
\frac{\Q^\leftarrow \circ H^\leftarrow (sx) -\Q^\leftarrow \circ H^\leftarrow (s)}{a\circ H^\leftarrow (s)}\cr
  &&\cr
&=&
\lim_{s\to\infty} 
\frac{
\Q^\leftarrow \Big( 
H^\leftarrow (s)+
\Big\{
\frac
{H^\leftarrow (sx) -H^\leftarrow (s)}{\sqrt{H^\leftarrow (s)}}
\Bigr\}\sqrt{H^\leftarrow(s)}\Bigr)
-\Q^\leftarrow \circ H^\leftarrow 
  (s)}{a\circ H^\leftarrow (s)}\cr
    &&\cr
&=&
\lim_{t\to\infty} 
\frac{
\Q^\leftarrow (t+\log x \cdot \sqrt t)
-\Q^\leftarrow (t)}{a(t)}\cr
  &&\cr
&=& 
 h^\leftarrow (-\log x).
\eean
 Thus, for $x>0$, using the form of $h^\leftarrow $
in \eqref{e:hinv},
\beqq\label{e:limV}
\lim_{s\to\infty} 
\frac{V (sx) -V 
  (s)}{a\circ H^\leftarrow (s)}
=\begin{cases}
-\frac 12 \log x, & \text{ if }\gamma=0,\\[5pt]
-\frac12 
\displaystyle{\bigg(\frac{x^{\gamma/2} -1}{\gamma/2}\bigg)}, & \text{ if } \gamma \neq 0.\end{cases}
\eeqq
From \cite[Theorem B.2.1, p.372]{dehaan:ferreira:2006}, we get
$a\circ H^\leftarrow (s) \in RV_{\gamma/2}$. 
Then multiply the limit relation in \eqref{e:limV} by $-1$ 
to see that the non-decreasing
function $-V$ is extended regularly varying at $\infty$ 
(\cite[p.127ff, p.139]{dehaan:ferreira:2006}).

We summarise the working up to \eqref{e:limV} as follows.

\begin{proposition}\label{p2} 
Assume  \eqref{e:Qinv}.
\begin{enumerate}[\rm (i)]
\item\label{item:1} When $\gamma<0$:
\beqq\label{e:item1}
V(x)=-V(\infty)- (-V(x)) \sim
a\circ H^\leftarrow (x)/|\gamma| \in RV_{-|\gamma|/2},\  {\rm as}\ x\to \infty;
\eeqq

\item When $\gamma=0$:  $-V \in {\bf \Pi}$ (or, equivalently,  $V \in {\bf   \Pi}_-$, \cite{dehaan:resnick:1979b})
with slowly varying auxiliary function $\frac 12 a\circ   H^\leftarrow
$.
\end{enumerate}
\end{proposition}

\begin{remark}\label{r3}
{\rm 
(i)\ 
Note $V$ being regularly varying with negative index in \eqref{e:item1} is consistent with $V$ being non-increasing.

(ii)\ 
Note that $\gamma >0$ cannot obtain in \eqref{e:limV}. 
The numerator on the left side of the limit is a
difference of two decreasing functions which as functions of $s$
approach $0$. The denominator is regularly varying with {\it
  positive\/} index and hence asymptotically increasing. So we cannot
get a non-trivial limit. See Remark \ref{rem:noPosGamma}.
}
\end{remark}

\subsection{Refining the centering and scaling for $\Y$.}\label{sec:refine}
Now we apply the material from Subsection \ref{subsec:implic} to
refine
the centering and scaling for $\Y$.
Recall that \eqref{e:LTcondit} or equivalently \eqref{e:Qinv} is in force throughout.
Depending on the range of $\gamma$, we may now simplify the form of the limit law for $\Y$ as follows.

\begin{proposition}\label{p3}
Suppose \eqref{e:LTcondit} holds with $\lchinv(y)$  as in \eqref{e:hinv}. Let $N_{\Gamma}$ be a standard normal random variable.
\begin{enumerate}[\rm (i)]
\item When $\gamma <0$: we may take $b_r  = \pibarinv(r)$, and then
\beqq\label{e:Del<0}
\frac{\Y}{b_r} \Rightarrow  e^{-N_\Gamma |\gamma|/2},\
{\rm as}\ r \to \infty.
\eeqq

\item When $\gamma=0$: we may take $b_r  = \pibarinv(r)$ and 
 $a_r = 2(\pibarinv(r-\sqrt{r}) - \pibarinv(r))$, and then
\beqq\label{e:0}
\frac{\Y - b_r}{a_r} \Rightarrow \frac{N_\Gamma}{2},\
{\rm as}\ r \to \infty. 
\eeqq
Furthermore
\be\label{e:ratioTo1}
a_r = o(b_r) \text{ and } \frac{\Delta^{(r)} }{b_r} \Rightarrow
1,\
{\rm as}\ r \to \infty.
\ee
\end{enumerate}
\end{proposition}

\medskip\noindent{\bf Proof of Proposition \ref{p3}:}\
(i)\ Take $\gamma <0.$
From \eqref{e:item1}, $V(x)\sim a\circ H^\leftarrow (x)/|\gamma|$, so 
\beqq\label{e:ba<0}
b_r=\Q^\leftarrow(r) =V(H(r)) \sim a_r/|\gamma|.\eeqq
Thus \eqref{e:YAN} can be written
$$\frac{\Y -b_r}{|\gamma|b_r} \Rightarrow h^\leftarrow (N_\Gamma),
$$
and hence, using \eqref{e:hinv},
$$\frac{\Y }{b_r} \Rightarrow 1+|\gamma|h^\leftarrow (N_\Gamma)
\eqdr e^{-|\gamma| N_\Gamma /2},
$$
which gives \eqref{e:Del<0}.

(ii) \ Take $\gamma=0$. 
From \eqref{e:hinv} with $\gamma=0$ and \eqref{e:Qinv} with $y=1$ we
get
$$\frac{2(\Q^\leftarrow (r-\sqrt r)-\Q^\leftarrow (r))}{a_r} \to 1,$$
and the choice of $a_r$ for \eqref{e:0} follows from  the convergence to types theorem. Since $V\in \bf{\Pi}_-$ with auxiliary function $a\circ
H^\leftarrow$ and the ratio of a non-negative  $\bf{\Pi}$ function  to its auxiliary
function tends to $\infty$, we have
$$\lim_{r\to\infty} \frac{b_r}{a_r}=\lim_{r\to\infty}
\frac{\Q^\leftarrow(r)}{a_r}
=\lim_{r\to\infty} \frac{V(H((r))}{a\circ H^\leftarrow\circ H(r)} =\infty.$$
Finally, dividing \eqref{e:0} by  $b_r/a_r$, which tends
to  $\infty$ as $r\to\infty$, yields a limit of $0$ which is tantamount to saying
$\Y/b_r\Rightarrow 1$.
\halmos

\begin{eg}\label{eg:stable}[Stable Subordinator]
 {\rm
 
To fix ideas, consider the case of the
  stable subordinator, where
\be\label{st}
\Q (x)=x^{-\alpha},\ x>0,\ 0<\alpha<1,\qquad \Q^\leftarrow
(y)=y^{-1/\alpha},\ y>0.
\ee
The numerator of the left side of \eqref{e:Qinv} is then, for $y\in \R$,
\begin{align*}
(r-y\sqrt r)^{-1/\alpha} -  r^{-1/\alpha}
&=
r^{-1/\alpha}\Bigl(\bigl(1-\frac{y}{\sqrt r}\bigr)^{-1/\alpha} -1\Bigr)\\
& \sim  r^{-1/\alpha -1/2} y/\alpha,\ {\rm as}\ r\to\infty,
\end{align*}
so  \eqref{e:Qinv} holds if we  take 
$$
b_r=r^{-1/\alpha},\quad 
a_r= 2r^{- 1/\alpha  - 1/2}/\alpha, 
\quad {\rm and}\quad
h^\leftarrow (y)=y/2.
  $$
  Thus we are in the $\gamma=0$ case.
  
Furthermore, recalling that $H^\leftarrow (x)=\frac 14 \log^2 x$, we get
$$
V(x)=\Q^\leftarrow \circ H^\leftarrow (x)=\Bigl( \frac 14 \log^2 x
\Bigr)^{-1/\alpha} \in  {\bf \Pi_-}.
$$
The auxiliary function corresponding to ${\bf \Pi}_-$-varying $V$ is
$$
a\circ H^\leftarrow (x) =\frac 2\alpha \Bigl( \frac 14 \log^2
x\Bigr)^{-1/\alpha -1/2},
$$ 
which is slowly varying at $\infty$ (as it should be), and for $x>0$
$$
\lim_{s\to\infty} \frac{V(sx)-V(s)  }{a\circ H^\leftarrow (s)   }
=-\frac 12 \log x.
$$
}\end{eg}
This completes the line-up of results needed for our analysis of $\Y$. Next we turn to the results needed for $\trimX$. 

\section{Further implications of the variation of  $\Q^\leftarrow$}\label{sec:more} 
Here we derive additional properties of $\Q^\leftarrow$ depending on
whether $\gamma <0$ or $\gamma=0$. These properties will be
 needed to replace random
centerings and scalings for $\trimX$ by deterministic normalizations in the following sections.

\subsection{Case $\gamma < 0$.} 
Suppose throughout that \eqref{e:LTcondit} holds with $h(x) = h_{\gamma}(x)$ for $\gamma < 0$ as in \eqref{e:h},  so by \eqref{e:Del<0} and \eqref{e:ba<0}  we can take  $a_r = |\gamma| b_r$ and
$b_r = \pibarinv(r)$, and have $\Y /b_r\Rightarrow e^{-N_\Gamma |\gamma|/2}.$
 Recall the   distribution  function  $G$ defined as $G(x)=e^{-\sqrt{\Q(x)}}$, $x>0$.
Then the following hold.

\begin{proposition}\label{prop:g<0} 
Assume  \eqref{e:LTcondit} holds with $\gamma<0$.
\begin{enumerate}[\rm (i)]
\item \label{prop:1}  For $p\geq 1$,
\beqq \label{e:momExpress}
\int_0^{b_r} u^p \nuu (\rmd u)  \sim \frac{2}{p|\gamma|} b_{r}^p \sqrt r, \ \text{as} \ r \to \infty.
\eeqq

In particular, when $p =2$,   $$\sigma^2(b_r) \sim
\frac{1}{|\gamma|}b_r^2 \sqrt{r}, \ \text{as} \ r \to \infty. $$

\item\label{item:sv}  $\Q (x)$ is slowly varying at $0$, 
  $\sigma^2(x)$ is regularly varying at $0$ with index 2
  and  $G(x)$ is regularly varying at 0   with index $1/|\gamma|$.
\end{enumerate}
\end{proposition}

\medskip\noindent{\bf Proof of Proposition \ref{prop:g<0}:}\
(i)\ 
Assume  \eqref{e:LTcondit} holds and keep $\gamma<0$ throughout.
To see (i), use  $\Q^\leftarrow =V\circ H$ from \eqref{e:defV}, where
$V$ is regularly varying at $\infty$ with index
$\gamma/2$ and $H$ is a ${\bf \Gamma}$ function with auxiliary
function $f(t)=\sqrt t$. Such a composition is again in the class
${\bf \Gamma}$ (\cite{resnickbook:2008,dehaan:ferreira:2006,dehaan:1970,dehaan:1974,bingham:goldie:teugels:1989,
resnickbook:2007}), 
so for $z\in \mathbb{R}$
$$
\frac{\Q^\leftarrow (r+\sqrt r z)}{\Q^\leftarrow (r)}=
\frac{V\left(\displaystyle{\frac{ H(r+\sqrt r z)}{H(r)}}H(r)\right)}{V(H (r))}
 \to
e^{z\gamma /2},
$$
or, equivalently, after a change of variable $w=-z|\gamma|/2$,
\beqq\label{e:Gvar}
\frac{\Q^\leftarrow (r+ \frac{2\sqrt r}{|\gamma|} w )}{\Q^\leftarrow (r)}  \to
e^{-w},\quad w\in \mathbb{R}.
\eeqq
The limit relation \eqref{e:Gvar} identifies the auxiliary function
of the non-increasing ${\bf \Gamma}$-varying
function 
$\Q^\leftarrow (x)$
as
 $f_1(r)=\frac{2}{|\gamma|}\sqrt r$. The function ${\bf \Gamma}$ class already
appeared in \eqref{e:limH} where we constructed a non-decreasing
function in ${\bf \Gamma}$. Likewise for any $p\geq 1$,
$(\Q^\leftarrow )^p \in {\bf \Gamma}$ with auxiliary function
$f_p(r)=\frac{2}{p|\gamma|} \sqrt r$. Auxiliary functions of
${\bf \Gamma}$-functions are unique up to asymptotic equivalence and 
also may be constructed in a canonical way
(see for example, \cite[page 19, eqn. 1.2.5]{dehaan:ferreira:2006},
\cite[p. 177, Corollary 3.10.5(b)]{bingham:goldie:teugels:1989}).
Therefore, we may  identify the auxiliary function of the
${\bf \Gamma}$-function $(\Q^\leftarrow )^p$  in two asymptotically equivalent ways:
\beqq \label{e:af}
f_p(r) \sim \frac{2}{p|\gamma|} \sqrt r 
\quad {\rm or}\quad 
f_p(r) \sim \frac{\int_r^\infty
  \bigl(\Q^\leftarrow (u)\bigr)^p 
  \rmd u}{\bigl(\Q^\leftarrow (r)\bigr)^p} ,\quad (r\to\infty).
\eeqq
Using the transformation theorem for integrals, we can write
(e.g. \cite[p. 301]{bremaud:1981})
$$\int_0^{\pibarinv(r)} u^p \nuu (\rmd u) =\int_r^\infty \bigl(\Q^\leftarrow (u)\bigr)^p
  \rmd u$$
and since $b_r=\Q^\leftarrow (r)$, applying \eqref{e:af} 
gives \eqref{e:momExpress}.

(ii)\
Invert the limit relation \eqref{e:Gvar} and change variables $s=b_r\to 0$ to get
$$
\lim_{s\to 0} \frac{\Q(sy)-\Q(s)}{(2/|\gamma|) \sqrt{\Q(s)}   }=-\log   y, \ y>0.
$$
Dividing by  $\Q(t)$ instead of $\sqrt{\Q(t)}$,   we get zero on the right side in the limit, which shows that $\Q(t) $ is slowly varying at $0$. Factoring as
\ben
\Q(sy)-\Q(s)=(\Q^{1/2}(sy)-\Q^{1/2}(s))(\Q^{1/2}(sy)+\Q^{1/2}(s))
\een
 and
using the slow variation of $\Q(x)$, hence of  $\Q^{1/2}(x)$, at $0$, gives the regular
variation of $e^{-\sqrt{\Q(x)}}$ at $0$ with index $1/|\gamma|$.\halmos

\subsection{Case $\gamma = 0.$}\label{sub:g=0}
Suppose \eqref{e:LTcondit} holds with $h(x) = h_{\gamma}(x)=2x$ for
$\gamma = 0$ as in \eqref{e:h}.
From   Proposition \ref{p3} we know in this case 
we may take $b_r  = \pibarinv(r)$ and  $a_r = 2(\pibarinv(r-\sqrt{r})
- b_r)$, and then $
(\Y - b_r)/a_r \Rightarrow {N_\Gamma}/{2}, $
where $N_{\Gamma}$ is a standard normal random variable. Also
$a_r/b_r \to 0$ and $\Delta^{(r)}/b_r\Rightarrow 1$. The following
proposition parallels Proposition \ref{prop:g<0} for the $\gamma=0$ case.
Recall the functions $H$ from \eqref{e:H} and $V=\Q^\leftarrow \circ
H^\leftarrow$ from \eqref{e:defV}, satisfying  $V^\leftarrow =H\circ \Q$ and $V \in {\bf \Pi_-} $ with slowly
varying auxiliary function $\frac 12 a\circ H^\leftarrow (s)$.

\begin{proposition}\label{p:=0}
Assume  that \eqref{e:LTcondit} holds with $\gamma=0$.
\begin{enumerate}[\rm (i)]
\item For $p\geq 1$, there exist ${\bf \Pi}$-varying functions $\pi_p(\cdot)$
  such that 
\beqq\label{e:mome}
\int_0^{b_r} u^p \Pi(du)=\pi_p(H(r))=\pi_p(e^{2\sqrt r})
\eeqq
where the slowly varying auxiliary function of $\pi_p$ is
$g_p(t)=\tfrac 12 V^p(t) \log t$.
\item As $r\to\infty$,
\beqq\label{e:truncVarNice}
\frac{\sigma^2(\Y)}{\sigma^2(b_r)} \Rightarrow 1.
\eeqq
\end{enumerate}
\end{proposition}

\medskip\noindent{\bf Proof of Proposition \ref{p:=0}:}\
 (i)\  For  $p\geq 1$ and $t>0$, recall $\hinv(y)= \tfrac{1}{4}\log^2 y$ and consider 
\begin{align}
\int_0^t u^p \nuu (\rmd u)=&\int_{\Q(t)}^\infty \bigl( \Q^\leftarrow
                         (s)\bigr)^p \rmd s  =\int_{\Q(t)}^\infty
                             \bigl( (V\circ H(s))^pds =\int_{H\circ
                             \Q(t)}^\infty  V^p(v) dH^\leftarrow (v)
\nonumber\\
=& \int_{V^\leftarrow (t) }^\infty V^p(v) \frac 12\log v \frac{\rmd v}{v}=
   \pi_p \bigl(V^\leftarrow (t)\bigr), \label{e:form} 
\end{align}
where we define
\be\label{e:pitilde}
 \pi_p (t)=\int_t^\infty V^p(v) \frac 12\log v \frac{\rmd v}{v}.
\ee
Now, $V$ is ${\bf \Pi}$-varying and hence slowly varying, so $V^p$ is
slowly varying, as is $\log v$. Thus the function $\pi_p (\cdot)$ is
the integral of a $-1$-varying function. The indefinite integral of a
$-1$-varying function is ${\bf \Pi}$-varying
(\cite{dehaan:1976, dehaan:ferreira:2006},
\cite[p. 30]{resnickbook:2008}). 
Thus $\pi_p \in {\bf \Pi}$ and the auxiliary function is $g_p(t)=\frac 12 V^p(t)\log t $. 

(ii)\ A ${\bf \Pi}$-varying function is always of larger order than its auxiliary function \cite[p.378]{dehaan:ferreira:2006}, so
\beqq\label{e:tooBig}
\lim_{t\to\infty}\frac{\pi_p(t)}{g_p(t)}=\infty.\eeqq

Now we apply these results with $p=2$.
Because of the representation in \eqref{e:form}, we invert the ${\bf
  \Pi}_-$-variation of $V(\cdot)$ in \eqref{e:limV} and get for $y>0,$
\beqq\label{e:limVinv}
\frac{V^\leftarrow \bigl(b_r-ya_r\bigr)}{V^\leftarrow (b_r)}
\to e^{2y},\ {\rm as}\ r\to\infty.
\eeqq
To  show \eqref{e:truncVarNice}, 
take the difference between numerator and denominator and use \eqref{e:form}:
\ben
\sigma^2 (\Y) -\sigma^2 (b_r)= \pi_2 \big( V^\leftarrow (\Y) \big)-  \pi_2
\big( V^\leftarrow (b_r) \big).
\een
From \eqref{e:YAN} write $(\Y-b_r)/a_r=\xi_r $ so that $\xi_r
	\Rightarrow N_\Gamma/2$ and remember $b_r=\Q^\leftarrow (r)$. The
	previous difference then becomes
	\ben
	   \pi_2 \big( V^\leftarrow (a_r\xi_r+b_r) \big)-  \pi_2  \big(
     V^\leftarrow (b_r) \big)\\
=  \pi_2 \Bigl(\frac{V^\leftarrow (a_r\xi_r+b_r)}{V^\leftarrow (b_r)}V^\leftarrow(b_r)\Bigr)-
 \pi_2 \bigl(V^\leftarrow (b_r)\bigr). 
\een
Applying the definition of ${\bf \Pi}$-variation and \eqref{e:limVinv}
we get
\beqq\label{e:justRight}
\frac{\sigma^2 (\Y) -\sigma^2 (b_r)}{g_2(V^\leftarrow (b_r))   }
\Rightarrow
N_\Gamma.\eeqq
Since 
$$\frac{\sigma^2(b_r)}{g_2(V^\leftarrow (b_r))}
=\frac{\pi_2(V^\leftarrow (b_r ))}{g_2(V^\leftarrow (b_r) )} \to
\infty,$$
by \eqref{e:tooBig}, we have proved \eqref{e:truncVarNice}, since if
we divide \eqref{e:justRight} by something of larger order (namely,
$\sigma^2(b_r)$), we get a limit of $0$. \halmos

 
 \section{Proofs of Theorems \ref{thm:trimXAN} and \ref{thm:joint} }\label{sect:proofs}
In this section we first prove the conditioned limit theorem,
Theorem \ref{thm:trimXAN},  using both random centering and scaling;
this is  followed by the  proof of Corollary \ref{cor:joint}; then we give the proof of Theorem  \ref{thm:joint}.
 The proof of Theorem   \ref{e:repmu} is deferred to Section \ref{p23}.

\medskip\noindent{\bf Proof of Theorem \ref{thm:trimXAN}:}\
Suppose $X$ is a driftless subordinator on
  $(0,\infty)$ with atomless  L\'evy measure $
\nuu(\cdot)$ on $(0,\infty)$ and its $r$th largest jump satisfies \eqref{Yconv} for some deterministic functions $a_r > 0$ and $b_r \in \R$.

Conditional on $\Y$,  we have that $\trimX$ is a subordinator whose
L\'evy measure is $\nuu_{|(0,\Y)}$, i.e., the measure $\Pi$ restricted to $(0,\Y)$  (e.g., \cite[Prop. 2.3, p.75]{resnick:1986b}).\footnote{Continuity of $\pibar$ is needed to apply Prop. 2.3 of 
\cite{resnick:1986b}; 
\cite{resnick:1986b} only gives the case $r=1$ but this is easily extended to $r\in\N$.}
So the conditional characteristic function (chf) of $\trimX $ is 
\ben
\EE \big( e^{\rmi \theta \,  \trimX} \big| \Y \big)=\exp \Big\{ 
\int_0^{\Y} \big( e^{\rmi \theta x}-1 \big)
  \nuu (\rmd x ) 
\Big\}, \ \theta \in \R,
\een
and the conditional chf of the centered and scaled $\trimX$ is 
\begin{align} \label{e:chf} 
&\EE\Big(\exp\Big\{ \rmi \theta \frac{\trimX -\mu(\Y)}{\sigma (\Y)} \Big\} \, \Big |\, \Y \Big) \nonumber \\
&=
\exp \Big\{ \int_0^{\Y} \Big(e^{\rmi \theta\frac{ u}{\sigma(\Y)}}-1- \rmi \theta \frac{u}{\sigma(\Y)}  \Big)  \nuu (\rmd u )  \Big\}.
\end{align}
Thus for \eqref{e:XCONV} it is enough to show  
\be\label{e:bdd0}
\Big|\int_0^{\Y} \Big(e^{\rmi \theta\frac{ u}{\sigma(\Y)}}-1- \rmi \theta \frac{u}{\sigma(\Y)}  \Big)  \nuu (\rmd u )   + \frac{1}{2}\theta^2 \Big| \to 0,\ \text{as } r \to \infty.
\ee

Noting that, by \eqref{e:defmusigma},
\[
 \int_0^{\Y} \frac{\theta^2}{2\sigma^2(\Y)}  u^2 \nuu (\rmd u) = \frac{1}{2}\theta^2,
\]
and using the inequality $|e^{\rmi \theta } - 1- \rmi \theta - \frac{(\rmi \theta)^2}{2}| \le |\theta|^3/3!$,
$\theta \in \R$, the lefthand side of $\eqref{e:bdd0}$ is seen to be 
\begin{align}\label{r7}
& \Big|\int_0^{\Y} \Big(e^{\rmi \theta\frac{ u}{\sigma(\Y)}}-1- \rmi \theta \frac{u}{\sigma(\Y)}  \Big)  \nuu (\rmd u )   - \int_0^{\Y} \Big(-\frac{1}{2}\theta^2 \Big) \frac{ u^2}{\sigma^2(\Y)}
\nuu (\rmd u)\Big| \nonumber \\
& \leq \frac{|\theta|^3}{3!} \int_0^{\Y}\frac{ u^3}{\sigma^3 (\Y)}\nuu
  (\rmd u)\cr
  & \leq \frac{|\theta|^3}{3!} \frac{\Y}{\sigma (\Y) }.
\end{align}

Next we show $\Y/ \sigma (\Y)$ converges to $0$ when  \eqref{Yconv}
holds. We
separate the analysis into cases according to whether the constant
$\gamma <0$ or $\gamma=0$ in \eqref{e:h}. 

\begin{enumerate}[(i)]
\item 
When $\gamma<0$, by \eqref{e:Del<0},  $\Y/b_r \Rightarrow Y : =
\exp(-N_\Gamma |\gamma|/2)$,  where $N_\Gamma$ is a 
standard normal random variable. Furthermore, $\sigma^2(t)$ is
regularly varying at $0$ with index $2$, so
\[
\frac{\sigma^2(\Y) }{\sigma^2 (b_r)} \to Y^2, \ \text{as} \ r \to \infty. 
\]
By \eqref{e:momExpress}, it is also true that $ \sigma^2(b_r)  \sim \frac{1}{|\gamma|}b_r^2 \sqrt r$. Then we have
\be\label{e:11}
 \frac{  (\Y)^2 }{\sigma^2(\Y) }
 \sim
 \frac{(\Y)^2}{b_r^2}  \cdot
\frac{b_r^2}{\sigma^2(b_r) }\cdot  \frac{\sigma^2(b_r)}{\sigma^2(\Y)}
=
 O_p(Y^2)\cdot \frac{|\gamma|}{\sqrt{r}}\cdot O_P(\frac{1}{Y^2})
 \Rightarrow  0, \ {\rm as}\ r\to\infty.
\ee

\item 
When $\gamma=0$, apply \eqref{e:ratioTo1} and then
Proposition \ref{p:=0}, and we have
$\Y /b_r\Rightarrow 1  $ and $\sigma^2(\Y) /\sigma^2 (b_r)
\Rightarrow 1.$ 
Therefore,  as in \eqref{e:11},
\begin{align*}
\frac{(\Y)^2 }{\sigma^2(\Y)} =&\frac{(\Y)^2}{b_r^2}  \cdot
\frac{b_r^2}{\sigma^2(b_r) }\cdot \frac{\sigma^2(b_r)}{\sigma^2(\Y)}
=(1+o_p(1)) \frac{b_r^2}{\sigma^2(b_r)} \Rightarrow 0, \ {\rm as}\ r\to\infty,
\end{align*}
where the  convergence to $0$ follows from \eqref{e:tooBig} for
the following reason: we can use \eqref{e:mome} (and recalling the definition of the function $g_p$ in  \eqref{e:mome}) to write
\begin{align*}
\lim_{r\to\infty} \frac{ \sigma^2 (b_r)  }{b^2_r
  }=&\lim_{r\to\infty} \frac{\pi_2 (H(r))}{(\Q^\leftarrow (r) )^2  }
      =\lim_{t\to\infty} \frac{\pi_2 (t)}{V^2(t)} 
=\lim_{t\to\infty} \frac{\pi_2 (t)}{V^2(t) \frac 12 \log t} \Bigl(\frac 12
      \log t\Bigr)\\
=& \lim_{t\to\infty} \frac{\pi_2 (t)}{g_2(t)} \Bigl(\frac 12 \log
   t\Bigr) =\infty.
\end{align*}
\end{enumerate}
Thus the righthand side of \eqref{r7} tends to 0, completing the proof of Theorem \ref{thm:trimXAN}.
\halmos

\medskip\noindent{\bf Proof of Corollary \ref{cor:joint}:}\
Define
$$\Phi_r^X=\frac{\trimX-\mu(\Y)}{\sigma(\Y)}, \quad 
\Phi_r^\Delta =\frac{\Y-b_r}{a_r}, $$
and suppose $f,g$ are non-negative continuous functions and bounded by 1. 
From  \eqref{e:onedim}  we have
$$
\EE g(\Phi_r^\Delta)  \EE f(N_X )  \to \EE(g(h^\leftarrow (N_\Gamma))    \EE f(N_X ).
$$
Then from \eqref{e:XCONV} and dominated convergence
$$
\EE f(\Phi_r^X ) g(\Phi_r^\Delta)= \EE \big\{ g(\Phi_r^\Delta) \EE^{\Y}\big( f(\Phi_r^X ) \big)
                            \big\}  \to  \EE \{f(N_X )\}\EE \big\{g \big( h^\leftarrow
                            (N_\Gamma) \big) \big\},$$
because
\begin{align*}
\bigl |\EE \big\{ g(\Phi_r^\Delta)& \EE^{\Y}(f(\Phi_r^X ) \bigr) \big\}     - \EE g(\Phi_r^\Delta) \EE f(N_X )   \bigr|\\
= &\, \Bigl | \EE \Bigl( g(\Phi_r^\Delta) 
\EE^{\Y} \big\{ f(\Phi_r^X )     -  \EE f(N_X ) \big\}  \Bigr)  \Bigr|\\
\leq &\,  \EE \big| \EE^{\Y}\big(f(\Phi_r^X )\big) -\EE f(N_X ) \big|\to 0, \ {\rm as}\ r\to\infty.
\end{align*}
This completes the proof of Corollary \ref{cor:joint}.
\halmos

\medskip\noindent{\bf Proof of Theorem \ref{thm:joint}:}\
(i)\  When $\gamma < 0$ 
we set $b_r = \pibarinv(r)$ and $a_r = |\gamma|b_r$, and then by \eqref{e:Del<0}
$$
\frac{\Y}{b_r} \Rightarrow Y : = e^{-N_\Gamma |\gamma|/2}, \ {\rm as}\ r\to\infty, 
$$
and the joint convergence in Corollary \ref{cor:joint} can be written
with a deterministic scaling via continuous mapping as
\begin{align*}
\Bigl(
\frac{\trimX-\mu(\Y)}{\sigma(b_r)}, \,
\frac{\Y}{b_r} \Bigr)
=&
\Bigl(
\frac{\trimX-\mu(\Y)}{\sigma(\Y)} \cdot \frac{\sigma(\Y)}{\sigma(b_r)},\,
\frac{\Y}{b_r} \Bigr)\\
 \Rightarrow &\bigl(N_X Y,\, Y\bigr)
=\bigl(N_X e^{-N_\Gamma |\gamma|/2},\, e^{-N_\Gamma |\gamma|/2}\bigr),
\end{align*}
where $(N_X,N_\Gamma)$ are iid standard normal random variables.

Now consider the effect of changing the random centering to a
deterministic one in the first component.
From \eqref{e:momExpress}, 
\be\label{r6}
\sigma(b_r)=\sqrt{\sigma^2(b_r)} \sim b_r r^{1/4}
\sqrt{\frac{1}{|\gamma|}}, \ {\rm as}\ r\to\infty.
\ee
Remember $a_r=|\gamma|b_r$ and convert \eqref{e:LTcondit} to vague
  convergence on $(0,\infty)$ to get
$$\frac{\nuu (b_r\rmd u)}{\sqrt r} \stackrel{v}{\to}  \frac{2}{|\gamma|}
\frac{\rmd u}{u},\quad u>0,
$$ 
and so
\begin{align*}
\frac{\mu (\Y) -\mu (b_r)}{b_r\sqrt r}
&=
\int_1^{\Y/b_r}
u\frac{\nuu(b_r\rmd u)}{\sqrt r}\\
&\Rightarrow 
\int_1^{Y} \frac{2}{|\gamma|}\rmd u=\frac{2}{|\gamma|}(Y-1).
\end{align*}
Since  $b_r\sqrt r/\sigma(b_r)\to\infty$ by \eqref{r6}, we have
\begin{align*}
\frac{\trimX-\mu(b_r)}{b_r\sqrt r}
&= \Bigl(
\frac{\trimX-\mu(\Y)}{b_r\sqrt r}\Bigr) +
\Bigl(\frac{\mu(\Y)-\mu(b_r)}{b_r\sqrt r} \Bigr)\\
&=
o_p(1) + \Bigl(\frac{\mu(\Y)-\mu(b_r)}{b_r\sqrt r} \Bigr)
\Rightarrow  \frac{2}{|\gamma|}  (Y-1), \ {\rm as}\ r\to\infty. 
\end{align*}

(ii)\
When $\gamma = 0$, by Proposition \ref{p:=0}, we have $\sigma(\Y) /
\sigma(b_r) \Rightarrow 1$ with $b_r = \pibarinv(r)$. Then from 
Corollary \ref{cor:joint}, 
\beqq\label{e:cheerylim}
\Bigl(
\frac{\trimX-\mu(\Y)}{\sigma(b_r)}, \,
\frac{\Y-b_r)}{a_r} \Bigr) \Rightarrow \bigl(N_X, \frac{N_\Gamma}{2}\bigr),
\eeqq
where $(N_X, N_\Gamma)$ are independent standard normal random
variables. By  Proposition \ref{p3} we may choose  $a_r = 2(\pibarinv(r-\sqrt{r}) - \pibarinv(r))$.
\halmos    

\section{Proof of Theorem   \ref{e:repmu}}\label{p23}
In this section we keep  $\gamma = 0$.  
 As displayed in \eqref{e:cheerylim}, we
may replace the random scaling for $\trimX$ by the deterministic
scaling $\sigma (b_r)$. 
We investigate what happens when we try to replace $\mu(\Y)$ with $\mu(b_r)$ in \eqref{e:cheerylim} by a  method similar to the one used in the proof of Proposition \ref{p:=0}.
As before we can  apply \eqref{e:form}, now with $p=1$, to  get
\ben
\mu (\Y) -\mu  (b_r)= \pi_1 \big( V^\leftarrow (\Y) \big)-  \pi_1
\big( V^\leftarrow (b_r) \big).
\een
Recall from \eqref{e:YAN} that we may write $(\Y-b_r)/a_r=\xi_r 
	\Rightarrow N_\Gamma/2$.  The 	previous difference thus becomes
	\ben
\pi_1 \big( V^\leftarrow (a_r\xi_r+b_r) \big)-  \pi_1  \big(
     V^\leftarrow (b_r) \big)=
      \pi_1 \Bigl(\frac{V^\leftarrow (a_r\xi_r+b_r)}{V^\leftarrow (b_r)}V^\leftarrow(b_r)\Bigr)-
 \pi_1 \bigl(V^\leftarrow (b_r)\bigr). 
\een
Applying the definition of ${\bf \Pi}$-variation and \eqref{e:limVinv}
we get
\beqq\label{e:ridMean}
\frac{\mu (\Y) -\mu (b_r)}{g_1(V^\leftarrow (b_r))   }
\Rightarrow
N_\Gamma.\eeqq

To replace $\mu(\Y)$ with $\mu(b_r)$  in \eqref{e:cheerylim} requires that the difference in 
\eqref{e:ridMean} be
compared with $\sigma(b_r)$. The cleanest result would be if the
difference were $o(\sigma (b_r))$ as $r\to\infty$, but this is not
always the case and   the final form of the joint limit with deterministic centering
and scaling in general depends on the behaviour of the limit
of 
\beqq\label{e:ratio}
\lim_{r\to\infty}\frac{\sigma^2(b_r)}{g_1^2(H(r))}
= \lim_{r\to\infty}\frac{ \pi_2(H(r))}{g_1^2(H(r))}
=\lim_{z\to\infty}\frac{\int_z^\infty \bigl( \Q^\leftarrow (v) \bigr)^2dv
}
{z \bigl( \Q^\leftarrow (z)  \bigr)^2 },
\eeqq
assuming there is indeed a limit.  Note that
 the ${\bf \Pi}$-function $\pi_2(\cdot) $ has  auxiliary function $g_2$
 and not $g_1^2$ so we cannot rely on \eqref{e:tooBig} here.

An easy example to show that
$(\mu(\Y)-\mu(b_r))/\sigma(b_r)$ does not always vanish 
is the stable subordinator from Example \ref{eg:stable}.
Recall from[\rm (i)] \eqref{st} we have, for $0<\alpha<1$,
$$\Q(x)=x^{-\alpha},\,x>0; \quad \Q^\leftarrow (v)=v^{-1/\alpha}, \, v>0.$$
The ratio on the right of \eqref{e:ratio} is in fact constant now:
$$
\frac{\int_z^\infty v^{-2/\alpha} \rmd v}{z(z^{-2/\alpha})}
=\frac{\alpha}{2-\alpha}.$$
More generally, if $\Q(x)=x^{-\alpha}L(x),\,x\downarrow 0$ is
regularly varying at $0$ with index $\alpha$, then $(\Q^\leftarrow
(z))^2 =z^{-2/\alpha}(L'(z))^2$, $z\to\infty$, for slowly varying
functions $L$ at $0$ and $L'$ at $\infty$. Then by Karamata's theorem
for integrals (eg. \cite[page 27]{bingham:goldie:teugels:1989})
$$
\lim_{z\to\infty}\frac{
\int_z^\infty \bigl( \Q^\leftarrow (v) \bigr)^2dv
}
{z \bigl( \Q^\leftarrow (z)
  \bigr)^2 }= \frac{\alpha}{2-\alpha}=c_\alpha.
  $$ 
Since $0<\alpha<1$, we have $0<c_\alpha<1$.
The converse half of Karamata's theorem
(\cite[p. 30]{bingham:goldie:teugels:1989})
tells us that if
\ben
\lim_{z\to\infty}\frac{\int_z^\infty\bigl(\Q^\leftarrow(v)\bigr)^2dv}
{ z\bigl(\Q^\leftarrow (z)  \bigr)^2} =c,\
{\rm for\ some}\ c\in (0,\infty),
  \een
then $ \bigl( \Q^\leftarrow (z)
  \bigr)^2 $ is regularly varying at $\infty$ with index 
  $-(c^{-1}+1)$ and $\Q^\leftarrow (z)$ is regularly varying with index
  $-(c^{-1}+1)/2$ at $\infty$.
Set $1/\alpha= (c^{-1}+1)/2$. Then for $\Q$ to correspond to a 
  subordinator, we need  $\alpha<1 $, which makes $c<1$.

Following this path leads to Theorem  \ref{e:repmu} as we state it in Section \ref{s2}, giving the joint limiting distribution of $\trimX$ and $\Y$ in this particular case. Based on the technology previously developed we can now prove that theorem.

\medskip\noindent{\bf Proof of Theorem   \ref{e:repmu}:}\
Suppose $\gamma=0$ and 
$\Q$ is regularly varying at $0$ with index  $-\alpha$,
$0<\alpha<1$. 
This happens iff $\Q^\leftarrow (z)$ is regularly varying at $\infty$
with index $-1/\alpha =-(1+c_\alpha^{-1})/2$,
where $c_\alpha= \alpha/(2-\alpha)$, 
or, equivalently, 
\beqq\label{e:assumption}
\lim_{z\to\infty}\frac{\int_z^\infty \bigl( \Q^\leftarrow (v) \bigr)^2dv
}{ \bigl( \Q^\leftarrow (z)  \bigr)^2 z                }
=\lim_{x\to0} \frac{ \int_0^xu^2 \nuu (\rmd u)}{x^2\Q(x)} =c_\alpha,\ {\rm where}\ c_\alpha \in (0,1).
\eeqq

Thus, suppose \eqref{e:assumption} holds. By \eqref{e:ratio}, we have $g_1(H(r))/\sigma(b_r)\to 1/\sqrt{c_\alpha}$. 
Then by \eqref{e:gamma0} and \eqref{e:ridMean},
\begin{align}
\frac{ \trimX -\mu(b_r)         }{\sigma(b_r)           }
&=
\frac{ \trimX -\mu(\Y)}{\sigma(b_r)} +\frac{ \mu(\Y) -\mu(b_r)}{\sigma(b_r)           }\label{e:divide}\\
&\eqdr 
N_X +o_p(1)+ \frac{ \mu(\Y) -\mu(b_r)}{g_1(H(r))} \cdot
   \frac{g_1 (H(r))   }{\sigma (b_r)} \nonumber\\
&\eqdr 
N_X+ N_\Gamma \cdot \frac{1}{\sqrt{c_\alpha}}+o_p(1).\nonumber\end{align}
Taking $r \to \infty$, this proves \eqref{e:depCase}.

When $c_\alpha=1$, the L\'evy measure property  that $$\int_0^1 u^2\nuu
(\rmd u)=\int_1^\infty \bigl( \Q^\leftarrow (s) \bigr)^2ds <\infty$$ 
(always) and the
usual version of Karamata's theorem, imply that  $\bigl(\Q^\leftarrow (x)\bigr)^2 $ is regularly
varying at infinity with index $-2$,  so $\Q^\leftarrow (x)$ is regularly
varying at infinity with index $-1$  and $\Q(x)$ is regularly varying at
0 with index $-1$. Conversely, if $\Q(x)$ is regularly varying at
0 with index $-1$, then \eqref{e:assumption} holds with $c_\alpha=1$ and so \eqref{e:depCase} holds with $c_\alpha=1$.

When \eqref{e:assumption} holds with $c_\alpha = 0$, then $\bigl(\Q^\leftarrow
(z) \bigr)^2 $ is rapidly varying at infinity  (\cite[p. 26]{dehaan:1970}) so the
same is true for $\Q^\leftarrow
(z)$, and, by inversion, $\Q(x) $ is slowly varying at $0$. The converse
holds as well: if $\Q$ is slowly varying at $0$ then
\eqref{e:assumption} holds with $c_\alpha=0$. Referring back to
\eqref{e:ratio} we find
$$\lim_{r\to\infty}\frac{\sigma^2(b_r)}{g_1^2(H(r))} =0.$$
Divide on the left side of \eqref{e:divide} by $g_1(H(r))$ instead of
$\sigma(b_r)$. Then by \eqref{e:ridMean} we see that 
 \eqref{e:depCase} becomes 
$$
\Bigl(
\frac{\trimX-\mu(b_r)}{g_1(H(r))      }, \,
\frac{\Y-b_r}{a_r} \Bigr) \Rightarrow 
\Bigl(N_\Gamma, \, \frac{N_\Gamma}{2}\Bigr),
$$
and unpacking the notation shows that $g_1\circ H(r)=b_r \sqrt r$ as claimed in Theorem  \ref{e:repmu}.\halmos

\section{Final thoughts}\label{s7}
One motivation for studying joint limit theorems as in Section \ref{s2} 
is to get information on limiting behaviour of ratios of a subordinator to its large jumps. See \cite{ipsen:maller:resnick:2017} and their references for related results and applications along these lines.

There are obvious open issues we leave for another day. Restricting
the investigation to subordinators clearly makes analysis easier but
we would like investigate what happens if we 
remove the assumption that the L\'evy process is
non-decreasing. This would presumably require analysis of the missing
case $\gamma >0$ which was necessarily absent from this paper.
We also would like to investigate functional weak limit theorems for
$(\trimX_t,\Y_t)$ as functions of $t$. Relevant to this, we note from 
\cite{buchmann:maller:resnick:2016}, Prop. 4.2, 
 that \eqref{e:LTcondit} implies, more generally,
\begin{equation}\label{e:onedimT}
\lim_{r\to\infty} \PP\Bigl( \frac{\Y_t-b_{r/t}}{a_{r/t}} \leq x \Bigr)
 = \PP\bigl(\Delta_t^{(\infty )} \leq x \bigr)=\Phi\bigl(\sqrt{t} h(x)\bigr),\ x\in\R,
\end{equation}
for each $t>0$. But $\trimX_t$ does not scale with $t$ in the same way 
as $\Y_t$, so generalisations of Theorems \ref{thm:trimXAN}--\ref{e:repmu}
are not straightforward in this respect.
As an incidental comment we note, though it's not mentioned in 
\cite{buchmann:maller:resnick:2016},  that \eqref{e:LTcondit} is in fact necessary and sufficient for \eqref{Yconv}.

\bibliographystyle{newapa}
\bibliography{bibfile.bib}

\end{document}